\UseRawInputEncoding 
\documentclass{article}

\usepackage{stmaryrd}
\usepackage{lmodern}
\usepackage{diagbox}
\usepackage{mathtools}
\usepackage{amssymb,amsmath}
\usepackage{graphicx}

\usepackage{mathrsfs}

\usepackage{color}
\usepackage{amsfonts}
\usepackage{epsfig}
\usepackage{graphics}

\usepackage{multicol}
\usepackage{lscape}
\usepackage{latexsym}
\def\proof{\noindent {\bf Proof }}
\def\Fq{{\mathbb F}_q}

\def\Zq-1{{\mathbb Z}_{q-1}}

\def\alp{{\alpha}}

\def\qed{~~\vrule height8pt width4pt depth0pt}

\def\Mcal{{\cal M}}
\def\Ical{{\cal I}}

\def\Pcal{{\cal P}}
\def\F2{{\mathbb F}_2}

\def\eps{{\varepsilon}}
\def\Mcal{{\cal M}}
\def\Ecal{{\cal E}}

\newtheorem{prop}{Proposition}

\newtheorem{thm}{Theorem}
\newtheorem{cor}{Corollary}

\newcommand{\fq}{{\mathbb F}_{q}}

\newcommand{\fqx}{{\mathbb F}_{q}^{*}}

\newcommand{\Scal}{{\cal S}}

\newcommand{\om}{{\omega}}

\begin{document}

\title {Enumeration of self-reciprocal irreducible monic polynomials with prescribed leading coefficients over a finite field }
 \author{
Zhicheng Gao\\
School of Mathematics and Statistics\\
Carleton University\\
Ottawa, Ontario\\
Canada K1S5B6\\
Email:~zgao@math.carleton.ca }
\maketitle

\begin{abstract}
A polynomial is called self-reciprocal (or palindromic) if the sequence of its coefficients is palindromic. In this paper we enumerate self-reciprocal irreducible  monic polynomials over a finite field with prescribed leading coefficients. Asymptotic expression with explicit error bound is derived, which implies that such polynomials of degree $2n$ always exist provided that the number of prescribed leading coefficients is slightly less than $n/4$. Exact expressions are also obtained for fields with two or three elements and with up to two prescribed leading coefficients.
\end{abstract}

\section{ Introduction}

In this paper we use our recent results from \cite{GKW21b} to enumerate self-reciprocal  irreducible monic polynomials  with prescribed leading coefficients. The following is a list of notations which will be used throughout the paper.

\begin{itemize}
\item $\fq$ denotes the finite field with $q$ elements, where $q=p^r$ for some prime $p$ and positive integer $r$.
\item $\Mcal_q$ denotes the set of monic polynomials over $\Fq$.
\item $\left[x^j\right]f(x)$ denotes the coefficient of $x^j$ in the polynomial $f(x)$.
\item For a polynomial $f$, $\deg(f)$ denotes the degree of $f$, and \\
$\displaystyle f^*(x)=x^{\deg(f)}f(1/x)$ is the {\em reciprocal} of $f$.
\item $\Pcal_q$ denotes the set of polynomials in $\Mcal_q$ with $f^*=f$. Polynomials in $\Pcal_q$ are called {\em self-reciprocal} or {\em palindromic}.
\item $\Ical_q\subseteq \Mcal_q$  denotes the set of irreducible monic polynomials.
\item $\Scal_q=\Ical_q\cap \Pcal_q$ denotes the set of  self-reciprocal  irreducible monic polynomials  over $\Fq$.
\item $\Scal_q(d)=\{f: f\in \Scal_q, \deg(f)=d\}$.
\item $\Ical_q(d)=\{f: f\in \Ical_q, \deg(f)=d\}$.
%\item $\Scal',\Scal''$ denote the set of polynomials in $\Scal$ with odd and even degrees, respectively.
\end{itemize}
Given non-negative integers $\ell,t$, and vectors $\vec{a}=(a_1,\ldots,a_{\ell})$ and $\vec{b}=(b_0,b_1,\ldots,b_{t-1})$, we also define
\begin{align*}
\Ical_q(d;\vec{a},\vec{b})&=\{f:f\in \Ical_q(d), \left[x^{d-j}\right]f(x)=a_{j}, 1\le j\le \ell, [x^j]f(x)=b_j, 0\le j\le t-1\},\\
\Scal_q(d;\vec{a})&=\{f:f\in \Scal_q(d), \left[x^{d-j}\right]f(x)=a_{j}, 1\le j\le \ell\}.
\end{align*}
Thus the vector $\vec{a}$ gives the $\ell$ {\em leading coefficients} of $f$, and $a_j$ is also called the $j$th {\em trace} of $f$. The vector $\vec{b}$ gives the $t$ {\em ending coefficients} of $f$ and $b_0$ is the {\em norm} of $f$.

There has been considerable interest in the study of $\Ical_q(d;\vec{a},\vec{b})$ and $\Scal_q(d;\vec{a})$, partly due to their applications in coding theory. Perhaps the most famous one is the Hansen-Mullen Conjecture \cite{HanMul92}, which says that irreducible polynomials with one coefficient prescribed (at a general position) always exist except for the two trivial forbidden cases. This conjecture was proved by Wan \cite{Wan97} for $d\ge 36$ or $q>19$ using character sums and Weil's bound. Similar questions were studied by Garefalakis and  Kapetanakis \cite{GarKap12} for self-reciprocal irreducible monic polynomials.

In this paper, our main focus is on the cardinality of $\Scal_q(d;\vec{a})$, which can be expressed in terms of $I_q(d;\vec{a},\vec{b}):=|\Ical_q(d;\vec{a},\vec{b})|$ as shown in the next section.

It is clear $I_q(d;\vec{a},\vec{b})=0$ when $d>1$ and $b_0=0$. Hence we assume $b_0\ne 0$ when $t\ge 1$. It is also easy to see that if $f(x)\in \Pcal_q$  has odd degree  then $f(-1)=0$, and consequently $\Scal_q(d)=\emptyset$ when $d>1$ is odd.  Thus we shall focus on self-reciprocal  irreducible monic polynomials of even degrees, and we set
\[
S_q(n;\vec{a}):=|\Scal_q(2n;\vec{a})|.
\]

The rest of the paper is organized as follows. In Section~\ref{main} we review some known results about  self-reciprocal polynomials and use it to derive a formula for $S_q(n;\vec{a})$ in terms of $I_q(d;\vec{a},\vec{b})$. We also state a formula for $I_q(d;\vec{a},\vec{b})$, which was obtained in the recent paper \cite{GKW21b}. These results will be used in Section~3 to derive asymptotic estimate for
$S_q(n;\vec{a})$.  In Sections~4 and 5 we obtain exact expressions for $S_q(n;\vec{a})$ with $q\in \{2,3\}$ and up to two prescribed leading coefficients.  Section~6 concludes the paper. Some tables of numerical values are provided in the last section as an appendix.

\section{Properties of self-reciprocal irreducible monic polynomials} \label{main}

The following results about self-reciprocal polynomials can be found in \cite{Mey90}.
\begin{prop}\label{prop:facts}
Let $f\in\Pcal_q$. Then
\begin{itemize}
\item[(a)] $f(x)$ has odd degree and is irreducible if and only if $f(x)=x+1$.
\item[(b)] $f$ has degree $2d$ if and only if there is a monic polynomial $g$ of degree $d$ with $g(0)=1$ such that
\begin{align}\label{eq:Quad}
    f(x)=x^dg\left(x+\frac{1}{x}\right).
\end{align}
      Moreover, when $d>1$, $g$ is irreducible if and only if either $f$ is irreducible or \begin{align}\label{eq:Quad1}
      f(x)=\frac{1}{h(0)}h(x)h^*(x)
      \end{align}
       for some $h\in \Ical_q(d)\setminus \Scal_q(d)$.
\end{itemize}
\end{prop}

Comparing the coefficients of both sides of  \eqref{eq:Quad}, we obtain
\begin{align}\label{eq:phi1}
 f_k=\sum_{j\le k/2}{d+2j-k\choose j}g_{k-2j},~~0\le k\le d.
\end{align}
We shall use $\phi_d:\fq^{\ell} \mapsto \fq^{\ell}$ to denote the mapping defined by \eqref{eq:phi1}
 from  $(g_1,\ldots,g_{\ell})$ to $(f_1,\ldots, f_{\ell})$.
 It is easy to see that $\phi_d$ is one-to-one as $g_k$ can be computed from $f_k$ recursively:
\begin{align}\label{eq:phiinverse}
 g_k=f_k-\sum_{0<j\le k/2}{d+2j-k\choose j}g_{k-2j}.
\end{align}

Equation~\eqref{eq:Quad1} corresponds to
\begin{align}\label{eq:psi}
 f_k=h_k+h_d^{-1}h_{d-k}+h_d^{-1}\sum_{j=1}^{k-1} h_jh_{d+j-k},\quad  h_d=h(0)\ne 0.
\end{align}
We note that the $\ell$ leading coefficients of $f(x)$ are determined by the $\ell$ leading coefficients and the $\ell+1$ ending coefficients of $h(x)$. Also for each given vector $\vec{b}:=(h_d,h_{d-1},\ldots,h_{d-\ell})$, system \eqref{eq:psi} gives a one-to-one mapping between $(f_1,\ldots,f_{\ell})$ and $(h_1,\ldots,h_{\ell})$.
We shall use $\psi_{\vec{b}}$ to denote this mapping, that is, $\psi_{(b_0,b_1,\ldots,b_{\ell})}$ is a bijection from $\fq^{\ell}$ to itself such that  $\psi_{(b_0,b_1,\ldots,b_{\ell})}(\vec{a})=\vec{c}$, where
\begin{align}
c_k=a_k+b_0^{-1}b_k+b_0^{-1}\sum_{j=1}^{k-1} a_jb_{k-j},~~1\le k\le \ell.\label{eq:psi1}
\end{align}

In the rest of the paper, we shall use the Iverson bracket $\llbracket P  \rrbracket$ which has value 1 if the predicate $P$ is true and 0 otherwise.

We now prove the following
\begin{thm} \label{thm:main} Suppose $n>1$. Let $\vec{c}=(c_1,\ldots,c_{\ell})\in \fq^{\ell}$  and $\vec{b}=(b_0,b_1,\ldots,b_{\ell})\in \fq^{\ell+1}$ with $b_0\ne 0$. Let $\phi_n$ and $\psi_{\vec{b}}$ be define above. Then
\begin{align}
S_{q}(n;\vec{c} )
&=\frac{1}{2}\sum_{\vec{a}\in \fq^{\ell} }\llbracket \psi_{(1,\vec{a})}(\vec{a})=\vec{c}\rrbracket S_q(n/2;\vec{a}) \nonumber\\
&~~~+I_q(n;\phi_n^{-1}(\vec{c}))-\frac{1}{2}\sum_{\vec{b}\in \fq^{\ell+1}}
I_q(n;\psi^{-1}_{\vec{b}}(\vec{c}),\vec{b}).\label{eq:main}
\end{align}
\end{thm}
\proof Proposition~\ref{prop:facts}(b) implies
\begin{align*}
S_{q}(n;\vec{c} )
&=I_q(n;\phi_n^{-1}(\vec{c}))\nonumber\\
&~~-\frac{1}{2}\left(\sum_{\vec{b}}
I_q(n;\psi^{-1}_{\vec{b}}(\vec{c}),\vec{b})-\sum_{\vec{a}\in \fq^{\ell} }\llbracket \psi_{(1,\vec{a})}(\vec{a})=\vec{c}\rrbracket S_q(n/2;\vec{a})\right),
\end{align*}
where $I_q(n;\phi^{-1}(\vec{c}))$ counts all polynomials formed by (1) and the second line corresponds to all polynomials formed by (2). \qed

Equation~\eqref{eq:main} immediately implies the following bounds
\begin{align}
 I_q(n;\phi_n^{-1}(\vec{c}))-\frac{1}{2}\sum_{\vec{b}}
I_q(n;\psi^{-1}_{\vec{b}}(\vec{c}),\vec{b})\le S_{q}(n;\vec{c})\le  I_q(n;\phi_n^{-1}(\vec{c})). \label{eq:Sbound0}
\end{align}

Theorem~\ref{thm:main} enables us to obtain expressions for $S_{q}(n;\vec{c})$ using known results about
$I_q(d;\vec{a},\vec{b})$.

For instance, if we take $\vec{c}$ to be the empty vector, we obtain the following formula for the total number of self-reciprocal irreducible monic polynomials of degree $2n$. This formula was first obtained by Carlitz \cite{Car67}; see also \cite{Coh69,Mey90}.
\begin{cor} We have
\begin{align}\label{eq:Sq2}
S_q(n)&=\left\{\begin{array}{ll}
\frac{1}{2n}\left(q^n+\llbracket 2\mid q\rrbracket-1\right),&\hbox{$n$ is a power of 2}, \\
\frac{1}{2n}\sum_{j\mid n}\llbracket2\nmid j\rrbracket\mu(j)q^{n/j},&\hbox{otherwise}.
\end{array}\right.
\end{align}
\end{cor}
\proof For self-completeness and illustration purpose, we include a short proof here.  A polynomial $f(x)=x^2+ax+1$ is reducible iff
$f(x)=(x+\alp)(x+1/\alp)$, that is, $a=\alp+1/\alp$ for some $\alp\in \fqx$. This implies
$S_q(1)=(q-1)/2$ when $q$ is odd, and $S_q(1)=q/2$ when $q$ is even. This gives \eqref{eq:Sq2} for $n=1$.

Now we assume $n\ge 2$ and write $n=2^ks$ with $k\ge 0$ and $s$ being an odd integer.
Taking $\vec{c}$ to be the empty vector in \eqref{eq:main} (that is, $\ell=0$), we obtain
\begin{align}\label{eq:L0}
S_{q}(n)&=\frac{1}{2}S_q(n/2)+\frac{1}{2}I_q(n),~~~n\ge 2,
\end{align}
which is \cite[(6.1)]{Coh69}.
Consequently,
\begin{align}\label{eq:SI1}
S_{q}(2^ks)&=2^{-k}S_q(s)+ \sum_{j=1}^{k} 2^{-j}I_q(n2^{1-j}).
\end{align}
Substituting the well-known formula (see, e.g.,\cite{LidNie97}; or apply our Theorem~\ref{thm:thm2} in the next section by taking $\ell=t=0$.)
\begin{align}
I_q(n)&=\frac{1}{n}\sum_{d\mid n}\mu(d)q^{n/d} \label{eq:Id}
\end{align}
into \eqref{eq:SI1}, and using
the property $\mu(2d)=-\mu(d)$ for any odd integer $d$, we obtain
\begin{align}
S_{q}(n)&=2^{-k}S_q(s)+ \sum_{j=1}^{k}\sum_{d\mid s2^{k+1-j}} \frac{1}{2n}I_q(n2^{1-j})\nonumber\\
&=\frac{s}{n}S_q(s)+\frac{1}{2n}\sum_{d\mid s}\sum_{j=1}^k\left(\mu(d)q^{n2^{1-j}/d}+\mu(2d)q^{n2^{-j}/d}\right)\nonumber\\
&=\frac{s}{n}S_q(s)+\frac{1}{2n}\sum_{d\mid s}\mu(d)\sum_{j=1}^k\left(q^{n2^{1-j}/d}-q^{n2^{-j}/d}\right)\nonumber\\
&=\frac{s}{n}S_q(s)+\frac{1}{2n}\sum_{d\mid s}\mu(d)\left(q^{n/d}-q^{s/d}\right).\label{eq:SI2}
\end{align}
When $s=1$, we have
\begin{align*}
S_{q}(n)&=\frac{1}{2n}\left(q+\llbracket 2\mid q\rrbracket-1\right)+\frac{1}{2n}\left(q^{n/d}-q\right),
\end{align*}
which gives the first line in \eqref{eq:Sq2}.

When $s>1$, we obtain from \eqref{eq:L0} and \eqref{eq:SI2}
\begin{align*}
S_q(n)&=\frac{s}{2n}I_q(s)+\frac{1}{2n}\sum_{d\mid s}\mu(d)\left(q^{n/d}-q^{s/d}\right).
\end{align*}
The second line in \eqref{eq:Sq2} follows immediately from \eqref{eq:Id}. \qed

We need a few more notations before stating the  formula from \cite{GKW21b} for $I_q(n;\vec{a},\vec{b})$.

We shall call two polynomials $f,g\in \Mcal_q$ {\em equivalent} with respect to $\ell,t$ if
\begin{align*}
\left[x^{\deg(f)-j}\right]f(x)&=\left[x^{\deg(g)-j}\right]g(x), 1\le j\le \ell, \\
\left[x^j\right]f(x)&=\left[x^j\right]g(x), 0\le j\le t-1.
\end{align*}
Thus polynomials are partitioned into equivalence classes. Let $\langle f\rangle$ denote the equivalence class represented by $f$.  It is known  \cite{GKW21b,Hay65,Hsu96} that the set $\Ecal^{\ell,t}$ of all equivalence classes forms an abelian group  under the multiplication
\[\langle f\rangle \langle g\rangle=\langle fg\rangle.
 \]
(When $t>0$, it is assumed that the constant term is nonzero.) It is also easy to see \cite{GKW21b,Hsu96} that
\[
\left|\Ecal^{\ell,t}\right|=(q-\llbracket t>0\rrbracket)q^{\ell+t-1}.
\]

 Since $\Ecal^{\ell,t}$ is abelian, it is isomorphic to a direct product of cyclic groups. Let $\xi_1,\ldots,\xi_f$ be the generators of these cyclic groups and denote their orders by $r_1,\ldots, r_f$, respectively. Thus each $\eps\in \Ecal^{\ell,t}$ can be written uniquely as
 \begin{align}
 \eps=\prod_{j=1}^f\xi_j^{e_j(\eps)}. \label{eq:exp}
 \end{align}

%We need a few more notations before stating the formula for $I_q(d;\eps)$, which is obtained recently in \cite{GKW21b}.

Let $\om_r=\exp(2\pi i/r)$ and $\eps\in \Ecal^{\ell,t}$.  Define
\begin{align}
\Ecal^{\ell,t}(d)&=\{\langle f \rangle: f\in \Mcal_q(d)\},\nonumber\\
%a(\eps,\eps')&=\prod_{h=1}^{f}  \om_{r_h}^{e_h(\eps)e_h(\eps')},\label{eq:a}\\
c(d;\eps) &= \sum_{\eps'\in \Ecal^{\ell,t}(d)}\prod_{j=1}^{f}  \om_{r_j}^{e_j(\eps)e_j(\eps')},\label{eq:cdj}\\
P(z;\eps)&=1+\sum_{d=1}^{\ell+t-1}c(d;\eps)z^d,\label{eq:Pj}\\
P(z)&=\prod_{\eps'\in \Ecal^{\ell,t} \setminus \{\langle 1\rangle\}} P(z;\eps'). \label{eq:Pz}
\end{align}

Following Granger's notation \cite{Gra19},  we set
\[
\rho_n(g):=\sum_{\rho} \rho^{-n},
\]
where the sum is over all the nonzero roots (with  multiplicity) of the polynomial $g\in {\mathbb C}[z]$.

Under the above notations, we can restate \cite[Theorem~3]{GKW21b} as follows.
\begin{thm}  \label{thm:thm2} Let $\Ecal$ denote $\Ecal^{\ell,t}$  and $\eps\in \Ecal$. We have
\begin{align}
I_q\left(n;\eps\right)&=\frac{1}{n}\sum_{k|n}\mu(k)\sum_{\delta\in \Ecal}\llbracket \delta^k=\eps\rrbracket
F_q(n/k;\delta),\label{eq:IF}
\end{align}
where
\begin{align}
F_q(n;\eps)&=\frac{q^n-\llbracket t>0 \rrbracket}{|\Ecal|}+\frac{n}{|\Ecal|}\sum_{\eps'\in \Ecal \setminus \{\langle 1\rangle\}}\prod_{j=1}^{f}  \om_{r_j}^{-e_j(\eps)e_j(\eps')}
[z^n]\ln  P(z;\eps')\nonumber\\
 &=\frac{q^n-\llbracket t>0 \rrbracket}{|\Ecal|}-\frac{1}{|\Ecal|}\sum_{\eps'\in \Ecal \setminus \{\langle 1\rangle\}}
\prod_{j=1}^{f}\om_{r_j}^{-e_j(\eps)e_j(\eps')}\rho_n(P(z;\eps')).\label{eq:Froot}
\end{align}
In particular, we have
\begin{align}
F_q(n;\langle 1\rangle)
 &= \frac{1}{|\Ecal|}\left(q^n-\llbracket t>0\rrbracket\right)-\frac{1}{|\Ecal|}
\rho_n(P(z)). \label{eq:trivial}
\end{align}
\end{thm}

\section{Asymptotic results}

In this section we use  Theorem~\ref{main} to prove the following bounds.
\begin{thm}\label{thm:thm3} Let $\vec{c}\in \fq^{\ell}$ and assume $1\le \ell\le n/2$. We have
\begin{align}
-\frac{\ell}{n}q^{\ell+1}q^{n/2}< S_q(n;\vec{c})-\frac{1}{2n}q^{n-\ell}
< \frac{\ell+1}{n}q^{\ell+1}q^{n/2}.\label{eq:Sbound}
\end{align}
Consequently, $S_q(n;\vec{c})>0$ whenever
\begin{align}\label{eq:L}
\ell\le \frac{n}{4}-\frac{\log_q (qn/2)}{2}.
\end{align}
\end{thm}
\proof Let $\vec{a}\in \fq^{\ell}$ and $\vec{b}\in \fq^{\ell+1}$ with $b_0\ne 0$. The following bounds follow immediately from \cite[Theoerm~2.1]{Coh05}):
\begin{align}
-\frac{\ell+1}{n} q^{n/2}&\le I_q(n;\vec{a})-\frac{q^{n-\ell}}{n}
\le \frac{\ell-1}{n} q^{n/2},\label{eq:Ibound0}\\
-\frac{2\ell+2}{n} q^{n/2}&\le I_q(n;\vec{a},\vec{b})-\frac{q^{-2\ell}}{n(q-1)}
\left(q^n-1\right)
\le \frac{2\ell}{n} q^{n/2}.\label{eq:Ibound1}
\end{align}
Using  \eqref{eq:main},  \eqref{eq:Ibound0} and \eqref{eq:Ibound1},   we obtain
\begin{align*}
S_q(n;\vec{c})&\le \frac{1}{n}\sum_{\vec{a}\in \fq^{\ell}}
\left(q^{n/2-\ell}+(\ell-1)q^{n/4}\right)+\frac{1}{n}q^{n-\ell}+\frac{\ell-1}{n}q^{n/2} \nonumber \\
&~~ -\frac{1}{2n}\sum_{\vec{b}\in \fq^{\ell+1},b_0\ne 0}
\left(\frac{q^{-2\ell}}{q-1}\left(q^n-1\right)-(2\ell+2)q^{n/2}    \right)\nonumber\\
&\le \frac{1}{2n}q^{n-\ell}+\frac{1}{n}\left(\ell-1+(\ell+1)(q-1)q^{\ell}\right)q^{n/2}
+\frac{\ell-1}{n}q^{\ell}q^{n/4}+\frac{q^{-\ell}}{2n}\nonumber\\
&< \frac{1}{2n}q^{n-\ell}+\frac{\ell+1}{n}q^{\ell+1}q^{n/2},
\end{align*}
which gives the upper bound in \eqref{eq:Sbound}.

For the lower bound, we use \eqref{eq:main},  \eqref{eq:Ibound0} and \eqref{eq:Ibound1} to obtain
\begin{align*}
S_q(n;\vec{c})&\ge \frac{1}{n}q^{n-\ell}-\frac{\ell+1}{n}q^{n/2} \nonumber \\
&~~ -\frac{1}{2n}\sum_{\vec{b}\in \fq^{\ell+1},b_0\ne 0}
\left(\frac{q^{-2\ell}}{q-1}\left(q^n-1\right)+2\ell q^{n/2}    \right)\nonumber\\
&> \frac{1}{2n}q^{n-\ell}-\frac{1}{n}\left(\ell+1+\ell(q-1)q^{\ell}\right)q^{n/2}\nonumber\\
&\ge \frac{1}{2n}q^{n-\ell}-\frac{\ell}{n}q^{\ell+1}q^{n/2}.
\end{align*}
It follows that $S_q(n;\vec{c})>0$ when
\[
q^{n/2}\ge 2\ell q^{2\ell+1}.
\]
So we may assume  $2\ell\le n/2$. Taking $\log_q$ on both sides, we complete the proof.
\qed

\medskip

 We note that Garefalakis and  Kapetanakis \cite{GarKap12} used character sums and Weil bound to derive estimate for the number of self-reciprocal irreducible monic polynomials with  one  coefficient prescribed at a general position.

\section{Exact results for $S_2(n;a)$ and $S_3(n;a)$}

In this section, we carry out more detailed calculations to obtain some explicit formulas for $S_2(n;a)$
and $S_3(n;a)$.

Substituting $\ell=1$ into \eqref{eq:phi1} and \eqref{eq:psi1}, we obtain
\begin{align*}
\phi_d(a)&=a,\\
\psi_{(b_0,b_1)}(a)&= a+b_0^{-1}b_1.
\end{align*}
It follows from Theorem~\ref{thm:main} that, for $n>1$,
\begin{align}
S_q(n;c)&=I_q(n;c)-\frac{1}{2}\sum_{b_1\in \fq,b_0\in \fqx} I_q(n;c-b_0^{-1}b_1,(b_0,b_1))\nonumber\\
&~~+\frac{1}{2}\sum_{a\in \fq}\llbracket 2a=c\rrbracket S_q(n/2;a)\nonumber\\
&=I_q(n;c)-\frac{1}{2}\sum_{b_1\in \fq,b_0\in \fqx} I_q(n;c-b_0^{-1}b_1,(b_0,b_1))\nonumber\\
&~~+\frac{1}{2}\llbracket 2\nmid q\rrbracket S_q(n/2;c/2)+\frac{1}{2}\llbracket 2\mid q, c=0\rrbracket S_q(n/2).\label{eq:S1}
\end{align}

The expression of $I_q(n;c)$ was first obtained by Carlitz \cite{Car52}. For self-completeness, we derive the result using Theorem~\ref{thm:thm2} again. In this case we note that $\ell=1,t=0$, each polynomial $P(z;\eps)$ is equal to 1, and hence the sum in \eqref{eq:Froot} is equal to 0. Thus
\[
F_q(n;\langle x+a\rangle)=q^{n-1}, ~~a\in \fq.
\]
Substituting this into \eqref{eq:IF}, we obtain
\begin{align}
I_q(n;c)&=\frac{1}{n}\sum_{j\mid n}\mu(j) q^{n/j-1}\sum_{a\in \fq}\llbracket \langle x+a\rangle^j=\langle x+c\rangle\rrbracket\nonumber\\
&=\frac{1}{qn}\sum_{j\mid n}\mu(j)q^{n/j}\sum_{a\in \fq}\llbracket ja=c\rrbracket.
\end{align}
Noting
\begin{align*}
\sum_{a\in \fq}\llbracket ja=0\rrbracket&=\llbracket p\nmid j\rrbracket+q\llbracket p\mid j\rrbracket
=1+(q-1)\llbracket p\mid j\rrbracket ,\\
\sum_{a\in \fq}\llbracket ja=c\rrbracket&=\llbracket p\nmid j\rrbracket,~~c\ne 0,
\end{align*}
 we obtain
\begin{align}
I_q(n;c)&=\frac{\llbracket c\ne 0\rrbracket}{qn}\sum_{j\mid n}\llbracket p\nmid j \rrbracket\mu(j)q^{n/j}\nonumber\\
&~~~+\frac{\llbracket c=0 \rrbracket}{qn}\sum_{j\mid n}(1+(q-1)\llbracket p\mid j \rrbracket)\mu(j)q^{n/j}.\label{eq:Iq1}
\end{align}

Next we derive a simple explicit expression for $S_2(n;1)$ and $S_2(n;0)$.

\begin{thm}\label{thm:S21} Let $\theta=\cos^{-1}\left(1/2\sqrt{2}\right)$. For $n\ge 1$, we have
\begin{align}
S_2(n;1)&=\frac{1}{4n}\sum_{j\mid n}\llbracket 2\nmid j\rrbracket\mu(j)
\left(2^{n/j}+1-(-1)^{n/j}2^{(n/2j)+1}\cos(n\theta/j)\right),\label{eq:S21s}\\
S_2(n;0)&=\frac{1}{4n}\sum_{j\mid n}\llbracket 2\nmid j\rrbracket\mu(j)
\left(2^{n/j}-1+(-1)^{n/j}2^{(n/2j)+1}\cos(n\theta/j)\right).\label{eq:S20s}
\end{align}
\end{thm}
\proof Since $x^2+x+1$ is irreducible over ${\mathbb F}_2$, we have $S_2(1;1)=1$ which agrees with the value obtained from \eqref{eq:S21s}.

For $n>1$, setting $q=2$ and $c=1$ in \eqref{eq:S1}, we obtain
\begin{align}
S_2(n;1)&=I_2(n;1)-\frac{1}{2}\left(I_2(n;1,(1,0))+I_2(n;0,(1,1)) \right).\label{eq:S21}
\end{align}

To find $I_2(n;a,(1,b))$, we first note that the corresponding group $\Ecal^{1,2}$ is generated by $\xi_1=\langle x+1\rangle$ and $\xi_2=\langle x^3+x+1\rangle$, and both $\xi_1$ and $\xi_2$ have order 2. We also have
\begin{align*}
\Ecal^{1,2}(1)=\{\xi_1\},\quad
\Ecal^{1,2}(2)=\{\langle x^2+1\rangle,\langle x^2+x+1\rangle\}=\{\langle 1\rangle,\xi_1\}.
\end{align*}
Using \eqref{eq:cdj} and \eqref{eq:Pj}, we obtain
\begin{align*}
c(1;\xi_1^{j_1}\xi_2^{j_2})&=(-1)^{j_1},\\
c(2;\xi_1^{j_1}\xi_2^{j_2})&=1+(-1)^{j_1},\\
P(z;\xi_1\xi_2^{j_2})&=1-z,\\
P(z;\xi_1^2\xi_2)&=1+z+2z^2\\
&=\left(1+\sqrt{2}e^{i\theta}z\right)\left(1+\sqrt{2}e^{-i\theta}z\right),
~~\theta=\cos^{-1}(1/2\sqrt{2}).
\end{align*}
It follows from \eqref{eq:Froot} that
\begin{align*}
F_2\left(n;\xi_1^{s_1}\xi_2^{s_2}\right)&=\frac{2^n-1}{4}-\frac{1}{4}\left((-1)^{s_1+s_2}+(-1)^{s_1}
+(-1)^{s_2+n}2^{n/2}\left(e^{in\theta}+e^{-in\theta}\right) \right)\nonumber\\
&=\frac{2^n-1}{4}-\frac{1}{4}\left((-1)^{s_1+s_2}+(-1)^{s_1}
+(-1)^{s_2+n}2^{(n/2)+1}\cos(n\theta)  \right).
\end{align*}
More explicitly, we have
\begin{align}
F_2(n;\langle 1\rangle)&=\frac{2^n-3}{4}-\frac{(-1)^n}{2}2^{n/2}\cos(n\theta),\label{eq:F200}\\
F_2(n;\xi_1)&=\frac{2^n+1}{4}-\frac{(-1)^n}{2}2^{n/2}\cos(n\theta),\label{eq:F210} \\
F_2(n;\xi_2)&= \frac{2^d-1}{4}+\frac{(-1)^n}{2}2^{n/2}\cos(n\theta),\label{eq:F201} \\
F_2(n;\xi_1\xi_2)&= \frac{2^n-1}{4}+\frac{(-1)^n}{2}2^{n/2}\cos(n\theta).\label{eq:F211}
\end{align}
It follows from \eqref{eq:IF} that
\begin{align}
I_2(n;\langle 1\rangle)&=\frac{1}{n}\sum_{j\mid n}\llbracket 2\nmid j\rrbracket\mu(j)F_2(n/j;\langle 1\rangle)\nonumber\\
&~~~+\frac{1}{n}\sum_{j\mid n}\llbracket 2\mid j\rrbracket\mu(j)\sum_{j_1,j_2\in \{0,1\}}F_2\left(n/j;\xi_1^{j_1}\xi_2^{j_2}\right)\nonumber\\
&=\frac{1}{4n}\sum_{j\mid n}\llbracket 2\nmid j\rrbracket\mu(j)\left(2^{n/j}-3-(-1)^{n/j}2^{(n/2j)+1}\cos(n\theta/j)\right)\nonumber\\
&~~~+\frac{1}{n}\sum_{j\mid n}\llbracket 2\mid j\rrbracket\mu(j)\left(2^{n/j}-1\right),\label{eq:I200}\\
I_2(n;\xi_1)&=\frac{1}{n}\sum_{j\mid n}\llbracket 2\nmid j\rrbracket\mu(j)F_2(n/j;\xi_1) \nonumber\\
&=\frac{1}{4n}\sum_{j\mid n}\llbracket 2\nmid j\rrbracket\mu(j)\left(2^{n/j}+1-(-1)^{d/j}2^{(n/2j)+1}\cos(n\theta/j)\right),\label{eq:I210}\\
I_2(n;\xi_2)&=I_2(n;\xi_1\xi_2)=\frac{1}{n}\sum_{j\mid n}\llbracket 2\nmid j\rrbracket\mu(j)F_2(n/j;\xi_2) \nonumber\\
&=\frac{1}{4n}\sum_{j\mid n}\llbracket 2\nmid j\rrbracket\mu(j)\left(2^{n/j}-1+(-1)^{n/j}2^{(n/2j)+1}\cos(n\theta/j)\right).\label{eq:I211}
\end{align}
Noting
\[
\xi_2=\langle x^3+x+1\rangle ,\quad  \xi_1\xi_2=\langle x^3+x^2+1\rangle,
\]
and substituting \eqref{eq:I211} into \eqref{eq:S21}, we obtain \eqref{eq:S21s}.

Using \eqref{eq:Sq2} and
\begin{align}
S_2(n;0)&=S_2(n)-S_2(n;1),
\end{align}
we obtain \eqref{eq:S20s}. \qed

\medskip

The expression for $S_3(n;a)$ is quite messy, and we will work out the expression of $I_3(n;a,(b_0,b_1))$ which can then be used to compute the numerical values of $S_3(n;a)$.

Substituting $q=3$ into \eqref{eq:S1}, (noting $2^{-1}=2$ in ${\mathbb F}_3$) we obtain, for $n>1$,
\begin{align}
S_3(n;0)&=\frac{1}{2}S_3(n/2;0)+I_3(n;0)\nonumber\\
&~~~ -\frac{1}{2}\left(I_3(n;0,(1,0))+I_3(n;0,(2,0)) \right)\nonumber\\
&~~~ -\frac{1}{2}\left(I_3(n;2,(1,1))+I_3(n;1,(2,1)) \right)\label{eq:S30}\\
&~~~ -\frac{1}{2}\left(I_3(n;1,(1,2))+I_3(n;2,(2,2)) \right).\nonumber
\end{align}

To obtain $I_3(n;a,(b_0,b_1))$, we first note that the corresponding group $\Ecal^{1,2}$ is generated by $\xi_1=\langle x+1\rangle$ and $\xi_2=\langle x^3+x+2\rangle$, which have orders 3 and 6, respectively. We also find
\begin{align*}
\Ecal^{1,2}(1)&=\{\langle x+1\rangle,\langle x+2\rangle\}=\{\xi_1,\xi_1^2\xi_2^3\},\\
\Ecal^{1,2}(2)&=\{\langle x^2+ax+b\rangle, 0\le a\le 2,1\le b\le 2\}=\{\langle 1\rangle,\xi_2^3,\xi_1,\xi_1\xi_2^5,\xi_1^2,\xi_1^2\xi_2\}.
\end{align*}

It follows from \eqref{eq:cdj} and \eqref{eq:Pj} that
\begin{align*}
c\left(1;\xi_1^{j_1}\xi_2^{j_2}\right)&=\om_3^{j_1}+\om_3^{2j_1}(-1)^{j_2},\\
c\left(2;\xi_1^{j_1}\xi_2^{j_2}\right)&=1+(-1)^{j_2}+\om_3^{j_1}+\om_3^{j_1}\om_6^{5j_2}+\om_3^{2j_1}+\om_3^{2j_1}\om_6^{j_2},\\
P\left(z;\xi_1^{j_1}\xi_2^{j_2}\right)&=1+c\left(1;\xi_1^{j_1}\xi_2^{j_2}\right)z+c\left(2;\xi_1^{j_1}\xi_2^{j_2}\right)z^2.
\end{align*}
If we use ${\bar g}(z)$ to denote the polynomial obtained from $g(z)\in {\mathbb C}[z]$ by taking the complex conjugate of the coefficients of $g(z)$, then we have
\begin{align*}
P\left(z;\xi_1^2\xi_2^{j_2}\right)={\bar P}\left(z;\xi_1\xi_2^{6-j_2}\right),\quad
P\left(z;\xi_1^3\xi_2^{j_2}\right)={\bar P}\left(z;\xi_1^3\xi_2^{6-j_2}\right).
\end{align*}
We also have
\begin{align*}
P\left(z;\xi_1\xi_2\right)&=1+i\sqrt{3}z,\\
P\left(z;\xi_1\xi_2^{2}\right)&=1-z+3z^2=\left(1-\frac{1+i\sqrt{11}}{2}z\right)\left(1-\frac{1-i\sqrt{11}}{2}z\right)\\
&=\left(1-\sqrt{3}e^{i\theta_1}z\right)\left(1-\sqrt{3}e^{-i\theta_1}z\right),
~~\theta_1=\cos^{-1}(1/2\sqrt{3}),\\
P\left(z;\xi_1\xi_2^3\right)&=1+i\sqrt{3}z,\\
P\left(z;\xi_1\xi_2^4\right)&=1-z,\\
P\left(z;\xi_1\xi_2^5\right)&=1+i\sqrt{3}z-3z^2=(1-i\sqrt{3}\om_3z)(1-i\sqrt{3}\om_3^2z),\\
P\left(z;\xi_1\xi_2^6\right)&=1-z   ,\\
P\left(z;\xi_1^3\xi_2\right)&=1+3z^2=(1+i\sqrt{3}z)(1-i\sqrt{3}z),\\
P\left(z;\xi_1^3\xi_2^2\right)&=1+ 2z+3z^2=\left(1+(1+i\sqrt{2})z\right)\left(1+(1-i\sqrt{2})z\right),\\
&=\left(1+\sqrt{3}e^{i\theta_2}z\right)\left(1+\sqrt{3}e^{-i\theta_2}z\right),
~~\theta_2=\cos^{-1}(1/\sqrt{3}),\\
P\left(z;\xi_1^3\xi_2^3\right)&=1,\\
P(z)&=\prod_{(j_1,j_2)\ne (3,6)}P\left(z;\xi_1^{j_1}\xi_2^{j_2}\right)\\
&=(1-z)^4(1+3z^2)^4(1-z+3z^2)^2(1+2z+3z^2)^2\\
&~~~\times(1-3z+3z^2)(1+3z+3z^2).
\end{align*}

Using \eqref{eq:Froot} and combining the conjugate pairs, we obtain
\begin{align}
F_3\left(n;\xi_1^{t_1}\xi_2^{t_2}\right)&=\frac{3^n-1}{18}
-\frac{1}{9}\left(\cos(2\pi(t_1+2t_2)/3)+\cos(2\pi t_1/3)\right)\nonumber\\
&-\frac{1}{9}3^{n/2}\cos((2t_1+t_2)\pi/3+n\pi/2)\nonumber\\
&-\frac{2}{9}3^{n/2}\cos(2(t_1+t_2)\pi/3)\cos(n\theta_1) \label{eq:F312f}\\
&-\frac{1}{9}3^{n/2}(-1)^{t_2}\cos(2t_1\pi/3+n\pi/2)\nonumber\\
&-\frac{1}{9}3^{n/2}\cos((2n-2t_1+t_2)\pi/3+n\pi/2)\nonumber\\
&-\frac{1}{9}3^{n/2}\cos((-2n-2t_1+t_2)\pi/3+n\pi/2)\nonumber\\
&-\frac{2}{9}3^{n/2}(-1)^{n/2}\llbracket 2\mid n\rrbracket\cos(t_2\pi/3)\nonumber\\
&-\frac{2}{9}3^{n/2}(-1)^n\cos(2t_2\pi/3)\cos(n\theta_2).\nonumber
\end{align}

Substituting $q=3$, $r_1=3$ and $r_2=6$ into \eqref{eq:IF}, and separating the six residue classes of $k$ (modulo 6), we obtain
\begin{align}
&~~I_3\left(n;\xi_1^{e_1}\xi_2^{e_2}\right)\nonumber\\
&=\frac{\llbracket e_1=e_2=0\rrbracket}{n} \sum_{k\mid n}\llbracket 6\mid k\rrbracket \mu(k)\left(3^{n/k}-1\right)\nonumber\\
&+\frac{1}{n}\sum_{k\mid n}\llbracket 6\mid k-1\rrbracket \mu(k)F_3\left(n/k;\xi_1^{e_1}\xi_2^{e_2}\right)\nonumber\\
&+\frac{1}{n}\sum_{k\mid n}\llbracket 6\mid k+1\rrbracket \mu(k)F_3\left(n/k;\xi_1^{-e_1}\xi_2^{-e_2}\right)\label{eq:I312f}\\
&+\frac{\llbracket 2\mid e_2\rrbracket}{n} \sum_{k\mid n}\llbracket 6\mid k-2\rrbracket \mu(k) \sum_{s_2\in \{0,3\}}F_3\left(n/k;\xi_1^{-e_1}\xi_2^{e_2/2+s_2}\right)\nonumber\\
&+\frac{\llbracket 2\mid e_2\rrbracket}{n} \sum_{k\mid d}\llbracket 6\mid k+2\rrbracket \mu(k) \sum_{s_2\in \{0,3\}}F_3\left(n/k;\xi_1^{e_1}\xi_2^{e_2/2+s_2}\right)\nonumber\\
&+\frac{\llbracket e_1=0,3\mid e_2\rrbracket}{n} \sum_{k\mid n}\llbracket 6\mid k+3\rrbracket \mu(k)\sum_{s_1,s_2\in \{0,1,2\}}F_3\left(n/k;\xi_1^{s_1}\xi_2^{e_2/3+2s_2}\right).\nonumber
\end{align}
We note that some cancellations in \eqref{eq:F312f} can be used to simplify the following sums:
\begin{align*}
&~~~\sum_{s_2\in \{0,3\}}F_3\left(n;\xi_1^{e_1}\xi_2^{e_2/2+s_2}\right)\\
&=\frac{3^{n}-1}{9}-\frac{2}{9}\left(\cos(2\pi (e_1+e_2)/3)+\cos(2\pi e_1/3) \right)\\
&~~~-\frac{4}{9}3^{n/2}\cos((2e_1+e_2)\pi /3)\cos(n\theta_1) \\
&~~~-\frac{4}{9}3^{n/2}(-1)^{n}\cos(e_2\pi/3)\cos(n\theta_2), \\
&~~~\sum_{s_1,s_2\in \{0,1,2\}}F_3\left(n;\xi_1^{s_1}\xi_2^{e_2/3+2s_2}\right)\\
&=\frac{3^{n}-1}{2}.
\end{align*}

In terms of the generators, we may rewrite \eqref{eq:S30} as
\begin{align}
S_3(n;0)&=\frac{1}{2}S_3(n/2;0)+I_3(n;0)-\frac{1}{2}\left(I_3(n;\langle 1\rangle)+I_3(n;\xi_2^3) \right)\nonumber\\
&~~~ -\frac{1}{2}\left(I_3(n;\xi_1^2\xi_2^4)+I_3(n;\xi_1\xi_2^5)+I_3(n;\xi_1\xi_2^2)+I_3(n;\xi_1^2\xi_2) \right)\label{eq:S30g},
\end{align}
which, together with \eqref{eq:F312f} and \eqref{eq:I312f}, gives a recursive way of computing $S_3(n;0)$.

Over ${\mathbb F}_3$, we have $2=-1$, and hence  $x\mapsto -x$ gives a bijection between self-reciprocal irreducible polynomials with trace 1 and those with trace $2$. Consequently
\begin{align}
S_3(n;1)&=S_3(n;2)=\frac{1}{2}(S_3(n)-S_3(n,0)). \label{eq:S31f}
\end{align}

Some numerical values of $I_3(n;a,(b_0,b_1))$ and $S_3(n;a)$ are calculated using \eqref{eq:F312f}--\eqref{eq:S31f}, which are given in Tables~1--4.

\section{Exact result for $S_2(n;a_1,a_2)$}

In this section, we derive some explicit expressions for computing  $S_2(n;a_1,a_2)$.

Substituting $\ell=2$ and $q=2$ into \eqref{eq:phi1} and \eqref{eq:psi}, we obtain
\begin{align*}
\phi_d(a_1,a_2)&=(a_1,a_2+d),\\
\psi_{(1,b_1,b_2)}(a_1,a_2)&=\left(a_1+b_1, a_2+b_2+b_1a_1\right).
\end{align*}
It follows from Theorem~\ref{thm:main} that, for $n>1$,
\begin{align}
S_2(n;c_1,c_2)&=\frac{1}{2}\sum_{a_1,a_2\in \{0,1\}}\llbracket c_1=0,a_1^2=c_2\rrbracket S_2(n/2;a_1,a_2)+I_2(n;c_1,n+c_2)\nonumber\\
&~~-\frac{1}{2}\sum_{b_1,b_2\in \{0,1\}} I_2(n;(c_1+b_1,c_2+b_2+b_1c_1+b_1),(1,b_1,b_2))\nonumber\\
&=\frac{\llbracket c_1=0\rrbracket }{2} S_2(n/2;c_2)+I_2(n;c_1,n+c_2)\label{eq:S22} \\
&~~-\frac{1}{2}\sum_{b_1,b_2\in \{0,1\}} I_2(n;(c_1+b_1,c_2+b_2+b_1c_1+b_1),(1,b_1,b_2)). \nonumber
\end{align}

An explicit expression for $I_2(n;a_1,a_2)$ was first obtained by Kuzmin \cite{Kuz90}. For self-completeness, we apply Theorem~\ref{thm:thm2} to obtain the following different looking  expressions.
\begin{prop}\label{ex:Q2L1} Let $\xi=\langle x+1\rangle$ be the generator of the group $\Ecal^{2,0}$. We have, for $0\le t\le 3$,
\begin{align}
&~~I_2\left(n;\xi^t\right)\nonumber\\
&=\frac{\llbracket t=0\rrbracket}{n} \sum_{k\mid n}\llbracket 4\mid k\rrbracket \mu(k)2^{n/k}
+\frac{\llbracket  2\mid t\rrbracket}{n} \sum_{k\mid n}\llbracket 4\mid k-2\rrbracket \mu(k)2^{(n/k)-1}\label{eq:I223f} \\
&+\frac{1}{n}\sum_{k\mid n}\llbracket 4\mid k-1\rrbracket \mu(k)\left(2^{(n/k)-2}-(-1)^{n/k}2^{(n/2k)-1}\cos\left(\frac{(n/k)-2t}{4}\pi\right)\right)\nonumber\\
&+\frac{1}{n}\sum_{k\mid n}\llbracket 4\mid k+1\rrbracket \mu(k)\left(2^{(n/k)-2}-(-1)^{n/k}2^{(n/2k)-1}\cos\left(\frac{(n/k)+2t}{4}\pi\right)\right).\nonumber
\end{align}
\end{prop}
\proof The group $\Ecal^{2,0}$ here is generated by $\xi=\langle x+1\rangle$ which has order 4. We have
\begin{align*}
\Ecal^{2,0}(1)&=\{\langle 1\rangle,\xi\},\\
c(1;\xi^j)&=1+i^j,\\
P(z;\xi^j)&=1+(1+i^j)z.
\end{align*}
Applying Theorem~\ref{thm:thm2}, we obtain
\begin{align}
F_2(n,\xi^t)&=2^{n-2}-\frac{1}{4}\sum_{j=1}^3i^{-jt}(-(1+i^j))^{n}\nonumber\\
&=2^{n-2}-\frac{1}{4}(-1)^n\left(i^{-t}(1+i)^n+i^{-3t}(1-i)^n \right)\nonumber  \\
&=2^{n-2}-\frac{1}{2}(-1)^n\Re\left(i^{-t}(1+i)^n\right)\nonumber\\
&=2^{n-2}-(-1)^{n}2^{(n/2)-1}\cos\left(\frac{(n-2t)\pi}{4}\right).
\end{align}

Substituting $q=2$, $r=4$ into \eqref{eq:IF},  and separating the four residue classes of $k$ (modulo 4), we  complete the proof. \qed

To compute $I_2(n;(a_1,a_2),(1,b_1,b_2))$, we need to find the generators of $\Ecal^{2,3}$. Using the computer algebra system {\em Maple}, we find that
$\Ecal^{2,3}$ is generated by $\xi_1=\langle x+1\rangle$ and $\xi_2=\langle x^5+x+1\rangle$. Both $\xi_1$ and $\xi_2$ have order 4.

We have
\begin{align*}
\Ecal^{2,3}(1)&=\{\xi_1\},\\
\Ecal^{2,3}(2)&=\{\langle x^2+ax+1\rangle: a\in \F2 \}=\{\xi_1^2,\xi_1^3\},\\
\Ecal^{2,3}(3)&=\{\langle x^3+ax^2+bx+1\rangle:a,b\in \F2\}=\{\langle 1\rangle,\xi_1\xi_2,\xi_1^2\xi_2^3,\xi_1^3\},\\
\Ecal^{2,3}(4)&=\{\langle x^4+ax^3+bx^2+cx+1\rangle: a,b,c\in \F2\}\\
&=\{\langle 1\rangle,\xi_1\xi_2^3,\xi_1,\xi_1^2,\xi_1^2\xi_2,\xi_1^3\xi_2^3,\xi_1^3,\xi_2\},
\end{align*}
and consequently,
\begin{align*}
c\left(1;\xi_1^{j_1}\xi_2^{j_2}\right)&=i^{j_1},\\
c\left(2;\xi_1^{j_1}\xi_2^{j_2}\right)&=(-1)^{j_1}+i^{3j_1},\\
c\left(3;\xi_1^{j_1}\xi_2^{j_2}\right)&=1+i^{j_1+j_2}+(-1)^{j_1}i^{3j_2}+i^{3j_1},\\
c\left(4;\xi_1^{j_1}\xi_2^{j_2}\right)&=1+i^{j_1+3j_2}+i^{j_1}+(-1)^{j_1}+(-1)^{j_1}i^{j_2}+i^{3j_1+3j_2}+i^{3j_1}+i^{j_2},\\
P\left(z;\xi_1^{j_1}\xi_2^{j_2}\right)&=1+c(1;\xi_1^{j_1}\xi_2^{j_2})z+c(2;\xi_1^{j_1}\xi_2^{j_2})z^2
+c(3;\xi_1^{j_1}\xi_2^{j_2})z^3+c(4;\xi_1^{j_1}\xi_2^{j_2})z^4,\\
P\left(z;\xi_1^{j_1}\xi_2^{j_2}\right)&={\bar P}\left(z;\xi_1^{4-j_1}\xi_2^{4-j_2}\right).
\end{align*}
Thus
\begin{align*}
P_{1,1}(z)&=P_{1,4}(z)=1+iz-(1+i)z^2,\\
P_{1,2}(z)&=P_{1,3}(z)=1+iz-(1+i)z^2+2(1-i)z^3 ,\\
P_{2,1}(z)&=1-z-2iz^3+4iz^4   ,\\
P_{2,2}(z)&=P_{2,4}(z)=1-z   ,\\
P_{4,1}(z)&=1+z+2z^2+2z^3+4z^4,\\
P_{4,2}(z)&=1+z+2z^2,\\
P(z)&=\prod_{(j_1,j_2)\ne (4,4)}P_{j_1,j_2}(z)\\
&=(1-z)^2(1-z^2-2z^3+2z^4)^2(1-z^2+2z^3-2z^4+8z^6)^2\\
&~~~\times(1+z+2z^2+2z^3+4z^4)^2\\
&~~~\times(1-2z+4z^6-16z^7+16z^8)(1+z+2z^2).
\end{align*}

Applying Theorem~\ref{thm:thm2} and combining conjugate pairs, we obtain
\begin{align}
F_2\left(n;\xi_1^{s_1}\xi_2^{s_2}\right)&=\frac{1}{16}\left(2^n-1\right)
-\frac{1}{16}\left((-1)^{s_1+s_2}+(-1)^{s_1}\right)  \nonumber \\
&~~~+\frac{n}{8}\Re\left(\left(i^{-s_1-s_2}+i^{-s_1}\right)[z^n]\ln(1+iz-(1+i)z^2)\right)\nonumber \\
&~~~+\frac{n}{8}\Re\left(\left(i^{-s_1-2s_2}+i^{-s_1-3s_2}\right)[z^n]\ln(1+iz-(1+i)z^2+2(1-i)z^3))\right)
\nonumber \\
&~~~+\frac{n}{8}\Re\left((-1)^{-s_1}i^{-s_2}[z^n]\ln(1-z-2iz^3+4iz^4)\right)\nonumber \\
&~~~+\frac{n}{8}\cos(\pi s_2/2)[z^n]\ln(1+z+2z^2+2z^3+4z^4) \nonumber\\
&~~~+\frac{n}{16}(-1)^{s_2}[z^n]\ln(1+z+2z^2).\label{eq:F2n}
\end{align}

Substituting $q=2$, $r_1=r_2=4$ into \eqref{eq:IF},  and separating the four residue classes of $k$ (modulo 4), we obtain
\begin{align}
&~~I_2\left(n;\xi_1^{e_1}\xi_2^{e_2}\right)\nonumber\\
&=\frac{\llbracket e_1=e_2=0\rrbracket}{n} \sum_{k\mid n}\llbracket 4\mid k\rrbracket \mu(k)\left(2^{n/k}-1\right)\nonumber\\
&+\frac{1}{n}\sum_{k\mid n}\llbracket 4\mid k-1\rrbracket \mu(k)F_2\left(n/k;\xi_1^{e_1}\xi_2^{e_2}\right)\nonumber\\
&+\frac{1}{n}\sum_{k\mid n}\llbracket 4\mid k+1\rrbracket \mu(k)F_2\left(n/k;\xi_1^{-e_1}\xi_2^{-e_2}\right)\label{eq:I223f}\\
&+\frac{\llbracket 2\mid e_1,2\mid e_2\rrbracket}{n} \sum_{k\mid n}\llbracket 4\mid k-2\rrbracket \mu(k) \sum_{s_1,s_2\in \{0,2\}}F_2\left(n/k;\xi_1^{e_1/2+s_1}\xi_2^{e_2/2+s_2}\right).\nonumber
\end{align}

As before, we may use some cancellations in \eqref{eq:F2n} to simplify the following sum:
\begin{align}
&~~\sum_{s_1,s_2\in \{0,2\}}F_2\left(n;\xi_1^{e_1/2+s_1}\xi_2^{e_2/2+s_2}\right)\nonumber\\
&=\frac{1}{4}\left(2^n-1
-(-1)^{(e_1+e_2)/2}-(-1)^{e_1/2}\right)\nonumber\\
&~~~+\frac{n}{2}\cos(\pi e_2/4)[z^n]\ln(1+z+2z^2+2z^3+4z^4) \nonumber\\
&~~~+\frac{n}{4}(-1)^{e_2/2}[z^n]\ln(1+z+2z^2) .
\end{align}

In terms of the generators $\xi$, $\xi_1$ and $\xi_2$, we may rewrite \eqref{eq:S22} as
\begin{align}
S_2(n;0,0)
&=\frac{1}{2} S_2(n/2;0)+I_2(n;0,n)\nonumber\\
&~~-\frac{1}{2}\sum_{b_1,b_2\in \{0,1\}} I_2(n;(b_1,b_1+b_2),(1,b_1,b_2))\nonumber\\
&=\frac{\llbracket 2\mid n\rrbracket}{2} S_2(n/2;0)+\llbracket 2\mid n\rrbracket I_2(n;\xi^0)+\llbracket 2\nmid n\rrbracket I_2(n;\xi^2)\nonumber\\
&~~-\frac{1}{2}\left(I_2(n;\xi_1^0)+I_2(n;\xi_1^2)+I_2(n;\xi_1^3\xi_2^2)+I_2(n;\xi_1\xi_2^2)\right), \label{eq:S200}\\
S_2(n;0,1)&=\frac{1}{2} S_2(n/2;1)+I_2(n;0,n+1)\nonumber\\
&~~-\frac{1}{2}\sum_{b_1,b_2\in \{0,1\}} I_2(n;(b_1,1+b_1+b_2),(1,b_1,b_2))\nonumber\\
&=\frac{\llbracket 2\mid n\rrbracket}{2} S_2(n/2;1)+\llbracket 2\mid n\rrbracket I_2(n;\xi^2)+\llbracket 2\nmid n\rrbracket I_2(n;\xi^0)\nonumber\\
&~~-\frac{1}{2}\left(I_2(n;\xi_2^2)+I_2(n;\xi_1^2\xi_2^2)+I_2(n;\xi_1)+I_2(n;\xi_1^3)\right), \label{eq:S201}\\
S_2(n;1,0)
&= I_2(n;1,n) \nonumber\\
&~~-\frac{1}{2}\sum_{b_1,b_2\in \{0,1\}} I_2(n;(1+b_1,b_2),(1,b_1,b_2))+\llbracket 2\nmid n\rrbracket I_2(n;\xi^2)\nonumber\\
&=\llbracket 2\mid n\rrbracket I_2(n;\xi)+\llbracket 2\nmid n\rrbracket I_2(n;\xi^3)\nonumber\\
 &~~-\frac{1}{2}\left(I_2(n;\xi_1\xi_2^3)+I_2(n;\xi_1^3\xi_2^3)+I_2(n;\xi_2)+I_2(n;\xi_1^2\xi_2)\right), \label{eq:S210} \\
 S_2(n;1,1)
&=I_2(n;1,n+1)\nonumber\\
&~~-\frac{1}{2}\sum_{b_1,b_2\in \{0,1\}} I_2(n;(1+b_1,1+b_2),(1,b_1,b_2))\nonumber\\
&=\llbracket 2\mid n\rrbracket I_2(n;\xi^3) +\llbracket 2\nmid n\rrbracket I_2(n;\xi)  \nonumber\\
&~~-\frac{1}{2}\left(I_2(n;\xi_1^3\xi_2)+I_2(n;\xi_1\xi_2)+I_2(n;\xi_1^2\xi_2^3)+I_2(n;\xi_2^3)\right).\label{eq:S211}
\end{align}

Some numerical values of $I_2(n;\xi_1^{e_1}\xi_2^{e_2})$ and $S_2(n,a_1,a_2)$ are computed using \eqref{eq:F2n}--\eqref{eq:S211}, which are given in tables~\ref{table5}--\ref{table7}.
For typographical convenience, the column headed by $\xi_1^{e_1}\xi_2^{e_2}$ in tables~\ref{table5} and \ref{table6} gives the values of $I_2(n;\xi_1^{e_1}\xi_2^{e_2})$.

\section{Conclusion} \label{conclusion}
We derived a general formula for the number $S_q(n;a_1,\ldots,a_{\ell})$ of self-reciprocal irreducible monic polynomials of degree $n$ over $\fq$ with leading coefficients $a_1,\ldots,a_{\ell}$. The general formula is expressed in terms of the number of irreducible monic polynomials with prescribed leading and ending coefficients. Using the general formula, we derived
exact expressions for $S_2(n;a)$, $S_3(n;a)$ and $S_2(n;a_1,a_2)$. Explicit error bounds were also obtained for $S_q(n;a_1,\ldots,a_{\ell})$ which imply that $S_q(n;a_1,\ldots,a_{\ell})>0$ when $\ell$ is slightly less than $n/4$.
In another paper \cite{Gao21}, we used Theorem~2 to improve the bounds for $I_q(n;\eps)$ given in \cite{Coh05,Hsu96}. Using the improved bounds for $I_q(n;\eps)$ and cancellations in the sums appeared in \eqref{eq:main}, we were able to significantly improve the bounds for $S_q(n;\eps)$.

\newpage
\section{Tables of numerical values}

\begin{table}[h]
\begin{center}
\begin{tabular}{|r|r|r|r|r|r|r|}
\hline
$n$&$I_3(n;\xi_2^0)$&$I_3(n;\xi_2)$&$I_3(n;\xi_2^2)$&$I_3(n;\xi_2^3)$&$I_3(n;\xi_2^4)$&$I_3(n;\xi_2^5)$\\
\hline
1 & 0& 0&0 &0&0 &0  \\
\hline
2 &1 &0 & 0& 0& 0& 0\\
\hline
3 & 0& 0 & 0& 0&1&1\\
\hline
4 & 0& 1 & 1& 2&1&1\\
\hline
5 &2& 3 & 3& 2&3&3\\
\hline
6&8& 8 & 6& 2&6&8\\
\hline
7 & 20&16 &16 & 20&16&16\\
\hline
8 &42 & 45& 44& 50&44&45 \\
\hline
9 &116&128 &128&116&119&119\\
\hline
10 & 334&328 &325 &320&325&328\\
\hline
11 &890&897&897&890&897&897\\
\hline
12 &2418&2447&2460&2504&2460&2447\\
\hline
13 &6848&6796&6796&6848&6796&6796\\
\hline
14 &18968&19002&18968&18884&18968&19002\\
\hline
15 &53072& 53103& 53103& 53072& 53249& 53249\\
\hline
16 &149370&149425&149380&149690&149380&149425\\
\hline
17 &422042&422019&422019&422042&422019&422019\\
\hline
18&1196484& 1195850& 1195362& 1195142& 1195362& 1195850\\
\hline
19 & 3398468& 3398404& 3398404& 3398468& 3398404& 3398404\\
\hline
20 &9682968& 9686128& 9685800& 9685264&9685800& 9686128\\
\hline
\end{tabular}
\end{center}
\caption{Values of $I_3(n;\xi_2^{s})$.}
\label{table1}
\end{table}
\begin{table}
\begin{center}
\begin{tabular}{|r|r|r|r|r|r|r|}
\hline
$n$&$I_3(n;\xi_1)$&$I_3(n;\xi_1\xi_2)$&$I_3(n;\xi_1\xi_2^2)$&$I_3(n;\xi_1\xi_2^3)$&$I_3(n;\xi_1\xi_2^4)$&$I_3(n;\xi_1\xi_2^5)$\\
\hline
1 & 1& 0&0 &0&0 &0  \\
\hline
2 &0 &0 & 0& 0& 0&1\\
\hline
3 & 0& 1&1& 0&0&1\\
\hline
4 & 1& 1 & 1& 2&1&0\\
\hline
5 &4& 3 & 2& 2&3&2\\
\hline
6&8& 8 & 6&6&6&7\\
\hline
7 &16&16 &18& 20&16&18\\
\hline
8 &41&45&50&46&44&44 \\
\hline
9 &120&119&121&120&128&121\\
\hline
10 & 337&328&310&332&325&328\\
\hline
11 &880&897&896&902&897&896\\
\hline
12 &2442&2447&2469&2448&2460&2476\\
\hline
13 &6856& 6796& 6816& 6800& 6796& 6816\\
\hline
14 &18974& 19002& 18902& 18940& 18986&19024\\
\hline
15 &53160& 53249& 53096& 53160& 53103& 53096\\
\hline
16 &149437& 149425& 149518& 149618& 149380& 149292\\
\hline
17 &422020& 422019& 422234& 421634& 422019& 422234\\
\hline
18&1195578& 1195850& 1195740& 1195698& 1195362& 1195861\\
\hline
19 &3398344& 3398404& 3398010& 3399380& 3398404& 3398010\\
\hline
20 &9684348& 9686128& 9685896& 9685668& 9685800& 9684248\\
\hline
\end{tabular}
\end{center}
\caption{Values of $I_3(n;\xi_1\xi_2^{s})$.}
\label{table2}
\end{table}

\begin{table}
\begin{center}
\begin{tabular}{|r|r|r|r|r|r|r|}
\hline
$n$&$I_3(n;\xi_1^2)$&$I_3(n;\xi_1^2\xi_2)$&$I_3(n;\xi_1^2\xi_2^2)$&$I_3(n;\xi_1^2\xi_2^3)$&$I_3(n;\xi_1^2\xi_2^4)$&$I_3(n;\xi_1^2\xi_2^5)$\\
\hline
1 & 0& 0&0 &1&0 &0  \\
\hline
2 &0 &1 & 0& 0& 0& 0\\
\hline
3 & 0& 1 &1& 0&1&0\\
\hline
4 & 1& 0 & 1& 2&1&1\\
\hline
5 &2&2 & 3&4&2&3\\
\hline
6&6& 7 & 6& 6&6&8\\
\hline
7 & 20&18 &16 & 16&18&16\\
\hline
8 &41 & 44& 44& 46&50&45 \\
\hline
9 &120& 121& 119& 120& 121& 128\\
\hline
10 & 337&328 &325 &332&310&328\\
\hline
11 &902&896&897&880&896&897\\
\hline
12 &2442&2476&2460&2448&2469&2447\\
\hline
13 &6800&6816&6796&6856&6816&6796\\
\hline
14 &18974&19024&18986&18940&18902&19002\\
\hline
15 &53160& 53096& 53249& 53160& 53096& 53103\\
\hline
16 &149437& 149292& 149380& 149618& 149518& 149425\\
\hline
17 &421634& 422234& 422019& 422020& 422234& 422019\\
\hline
18&1195578& 1195861& 1195362& 1195698& 1195740& 1195850\\
\hline
19 & 3399380& 3398010& 3398404& 3398344& 3398010& 3398404\\
\hline
20 &9684348& 9684248& 9685800& 9685668& 9685896& 9686128\\
\hline
\end{tabular}
\end{center}
\caption{Values of $I_3(n;\xi_1^2\xi_2^{s})$.}
\label{table3}
\end{table}

\begin{table}
\begin{center}
\begin{tabular}{|r|r|r|}
\hline
$n$&$S_3(n;0)$&$S_3(n;1)$\\
\hline
1 & 1& 0\\
\hline
2 &0 &1\\
\hline
3 & 0& 2\\
\hline
4 & 4& 3\\
\hline
5 &10&7\\
\hline
6&20&208\\
\hline
7 &48&54\\
\hline
8 &132&139 \\
\hline
9 &368&362\\
\hline
10 &1000&976\\
\hline
11 &2686&2683\\
\hline
12 &7340&7400\\
\hline
13 &20400& 20460\\
\hline
14 &57000& 56910\\
\hline
15 &159584& 159352\\
\hline
16 &448396& 448407\\
\hline
17 &1265650& 1266295\\
\hline
18&3586820& 3587420\\
\hline
19 &10196064& 10194882\\
\hline
20 &29058328& 29055640\\
\hline
\end{tabular}
\end{center}
\caption{Values of $S_3(n;0)$ and $S_3(n;1)=S_3(n;2)$.}
\label{table4}
\end{table}

\begin{table}
\begin{center}
\begin{tabular}{|r|r|r|r|r|r|r|r|r|}
\hline
$n$&$\xi_2^0$&$\xi_2$&$\xi_2^2$&$\xi_2^3$&$\xi_1$
&$\xi_1\xi_2$&$\xi_1\xi_2^2$&$\xi_1\xi_2^3$\\
\hline
1 & 0& 0& 0& 0& 1& 0& 0& 0\\
\hline
2 &0 &0& 0& 0& 0& 0&0& 0\\
\hline
3 &0 &0& 0& 0& 0& 1&0& 0\\
\hline
4  &0 &1& 0& 0& 0& 0&0& 1\\
\hline
5  &0 &0& 1& 0& 0& 0&1& 0\\
\hline
6&1 &1& 0& 0& 0& 1&1& 1\\
\hline
7 &2 &1& 1& 1& 2& 1&0& 1\\
\hline
8 &1 &2&2& 1& 2& 2&2& 2\\
\hline
9 &4 &4& 3& 4& 6& 2&3& 4\\
\hline
10 &5 &5&6& 8& 5& 7&7& 5\\
\hline
11 &10 &11& 13& 11& 12& 11&14& 11\\
\hline
12 &23 &24& 18& 20& 21& 20&20&24\\
\hline
13 &36 &42&35&42&36&36&41&42\\
\hline
14 &73&75&70&70&76&73&73&75\\
\hline
15 &138&137&138&137&134&133&137&137\\
\hline
16 &243&262&258&245&238&262&262&262\\
\hline
17 &484&488&475&488&482&486&479&488\\
\hline
18&930&889&894&913&913&912&912&889\\
\hline
19 &1722&1719&1725&1719&1728&1743&1722&1719\\
\hline
20 &3327 &3271&3234&3275&3260&3288&3288&3271\\
\hline
\end{tabular}
\end{center}
\caption{Values of $I_2(n;\xi_2^s)$ and $I_2(n;\xi_1\xi_2^s)$, $0\le s\le 3$.}
\label{table5}
\end{table}

\begin{table}
\begin{center}
\begin{tabular}{|r|r|r|r|r|r|r|r|r|}
\hline
$n$&$\xi_1^2$&$\xi_1^2\xi_2$&$\xi_1^2\xi_2^2$&$\xi_1^2\xi_2^3$&$\xi_1^3$
&$\xi_1^3\xi_2$&$\xi_1^3\xi_2^2$&$\xi_1^3\xi_2^3$\\
\hline
1 & 0& 0& 0& 0& 0& 0& 0& 0\\
\hline
2 &0 &0& 0& 0& 1& 0&0& 0\\
\hline
3 &0 &0& 0& 1& 0& 0&0& 0\\
\hline
4  &0 &0& 0& 0&1& 0&0& 0\\
\hline
5  &0 &1& 1& 0& 0& 0&1& 1\\
\hline
6&0 &1& 0&1& 0&0&1& 1\\
\hline
7 &0 &2& 1& 1& 2& 1&0&2\\
\hline
8 &2 &2&2&2&3&1&2& 2\\
\hline
9 &4 &4& 3&2&2&4&3& 4\\
\hline
10 &4 &7&6&7& 5&8&7&7\\
\hline
11 &12 &12& 13& 11& 8& 11&14& 12\\
\hline
12 &22 &20& 18& 20& 25& 20&20&20\\
\hline
13&48&41&35&36&36&42&41&41\\
\hline
14 &72&73&70&73&72&70&73&73\\
\hline
15 &136&134&138&133&142&137&137&134\\
\hline
16 &242&262&258&262&255&245&262&262\\
\hline
17 &492&467&475&486&486&488&479&467\\
\hline
18&908&912&894&912&917&913&912&912\\
\hline
19 &1716&1728&1725&1743&1716&1719&1722&1728\\
\hline
20 &3246 &3288&3234&3288&3256&3275&3288&3288\\
\hline
\end{tabular}
\end{center}
\caption{Values of $I_2(n;\xi_1^2\xi_2^s)$ and $I_2(n;\xi_1^3\xi_2^s)$, $0\le s\le 3$.}
\label{table6}
\end{table}

\begin{table}
\begin{center}
\begin{tabular}{|r|r|r|r|r|}
\hline
$n$&$S_2(n;0,0)$&$S_2(n;0,1)$&$S_2(n;1,0)$&$S_2(n;1,1)$\\
\hline
1 & 0& 0 &0 & 1\\
\hline
2 & 0& 0 & 0& 1\\
\hline
3 & 1& 0 & 0& 0\\
\hline
4 & 1& 0 & 0& 1\\
\hline
5 & 1& 0 & 1& 1\\
\hline
6& 1& 2 & 1& 1\\
\hline
7 & 3& 2 & 2& 2\\
\hline
8 &3&4 &4&5\\
\hline
9 & 6&8 &5&9\\
\hline
10 & 13&14 & 12&12\\
\hline
11 &23&22&22&26\\
\hline
12 &44&40 &41&45\\
\hline
13 &77&84 &77&77\\
\hline
14 & 145&146 & 149&145\\
\hline
15 &267&274 &279&271\\
\hline
16 &507&524 &500&517\\
\hline
17 &953&976 &965&961\\
\hline
18& 1802&1824& 1825&1829\\
\hline
19 &3471&3438 &3438& 3450\\
\hline
20 &6546&6576&6548&6544\\
\hline
\end{tabular}
\end{center}
\caption{Values of $S_2(n;0,0),S_2(n;0,1),S_2(n;1,0),S_2(n;1,1)$.}
\label{table7}
\end{table}

\end{document}